\documentclass[11pt]{amsart}
\usepackage{amssymb, amsmath, amsthm, amsfonts, bm}
\usepackage{verbatim}
\usepackage{bbm}
\usepackage{enumerate}
\usepackage{dsfont}
\usepackage[mathscr]{eucal}
\usepackage{color}
\usepackage{ulem}
\usepackage{cancel}
\usepackage{mathtools}

\usepackage{mathabx}

\def\beq{\begin{equation}}
\def\endeq{\end{equation}}
\def\lesim{\lesssim}

\def\vector{\overrightarrow}

\def\be#1{\begin{equation}\label{ #1}}
\def\endal{\end{align}}
\def\bas{\begin{align*}}
\def\eas{\end{align*}}
\def\bi{\begin{itemize}}
\def\ei{\end{itemize}}

\def\emph#1{{\it #1}}
\def\textbf#1{{\bf #1}}

\def\R{\mathbb{R}}

\def\Z{\mathbb{Z}}

\newcommand{\Be}{\begin{equation}}
\newcommand{\Ee}{\end{equation}}

\newcommand{\Bm}{\begin{multline}}
\newcommand{\Em}{\end{multline}}
\newcommand{\Bea}{\begin{eqnarray}}
\newcommand{\Eea}{\end{eqnarray}}
\newcommand{\Beas}{\begin{eqnarray*}}
\newcommand{\Eeas}{\end{eqnarray*}}
\newcommand{\Benu}{\begin{enumerate}}
\newcommand{\Eenu}{\end{enumerate}}
\newcommand{\Bi}{\begin{itemize}}
\newcommand{\Ei}{\end{itemize}}

\def\intslash{\rlap{\kern  .32em $\mspace {.5mu}\backslash$ }\int}
\def\qsl{{\rlap{\kern  .32em $\mspace {.5mu}\backslash$ }\int_{Q_x}}}

\def\emph#1{{\it #1 }}

\def\diag{{\text{\rm diag}}}
\def\supp{{\text{\rm supp}}}

\numberwithin{equation}{section}

\theoremstyle{plain}
\newtheorem{thm}{Theorem}[section]

\newtheorem*{thm*}{Theorem}
\newtheorem*{conj*}{Conjecture}

\newtheorem*{openproblem*}{Open Problem}

\begin{document}

\title[Fourier supports on moment curves]{$L^p$ integrability of functions with Fourier support on a smooth space curve}

\author{Shaoming Guo, Alex Iosevich, Ruixiang Zhang, \\
Pavel Zorin-Kranich}

\begin{abstract} We prove that if $f\in L^p(\R^k)$ with $p<(k^2+k+2)/2$ satisfies that $\widehat{f}$ is supported on a small perturbation of the moment curve in $\R^k$, then $f$ is identically zero.
This improves the more general result in \cite{AN04} and the exponents are sharp in all dimensions.
 In the process, we develop a mechanism that should lead to further progress on related problems.
\end{abstract}

\maketitle
\date{\today}

\section{Introduction}

A classical result due to Agranovsky and Narayanan (\cite{AN04}) says that if the Fourier transform of an integrable function is supported on a smooth compact $s$-dimensional sub-manifold of ${\Bbb R}^k$, $1 \leq s \leq k-1$, and $f \in L^{p}({\Bbb R}^k)$ for $p<\frac{2k}{s}$, then $f$ is identically $0$. When $s \ge \frac{k}{2}$, this result is essentially sharp (see \cite{AN04} for details). For example, if $s=k-1$ and the manifold is the unit sphere, a direct calculation involving the integrability of the Fourier transform of the surface measure on the sphere yields sharpness. When $s<\frac{k}{2}$, the situation changes drastically and the $\frac{2k}{s}$ threshold is never sharp. For example, if $s=1$, and the manifold is a compact piece of the moment curve $\{(t,t^2, \dots, t^{k}): t \in [0,1] \}$, a simple scaling argument (\cite{B20}) reduces the problem to the integrability results in \cite{BGGIST07}. The main purpose of this paper is to provide a robust method for the study of the problem under consideration, which should open the door to a variety of related results.

\vskip.125in 

\begin{thm}\label{main_thm}
Fix a positive integer $k\ge 2$ and a real number $1\le p< (k^2+k+2)/2$. If $f\in L^p(\R^k)$ satisfies that $\widehat{f}$ is supported on 
\begin{equation}\label{moment_curve}
\{\vector{\gamma}(t):=(t, t^2, \dots, t^k) \big| t\in [0, 1]\},
\end{equation}
then $f\equiv 0$. Moreover, the conclusion is still valid for small perturbations of this curve. 
\end{thm}

Arkhipov, Chubarikov and Karatsuba \cite{ACK87} proved that 
\begin{equation}
\int_0^1 e^{i(t x_1+t^2 x_2+\dots+t^k x_k)} dt,
\end{equation}
as a function of $(x_1, \dots, x_k)$, belongs to $L^p(\R^k)$ if and only if 
$$p>(k^2+k+2)/2, $$ which establishes the sharpness of Theorem \ref{main_thm} up to the endpoint.

The case $k=2$ in Theorem \ref{main_thm} was obtained by Agranovsky and Narayanan \cite{AN04}, as we noted above. We remark here that in the result of Agranovsky and Narayanan \cite{AN04}, there is no curvature assumptions on the underlying  $C^1$ manifold. For instance, their result applies to both hyperplanes and hypersurfaces like paraboloids. However, we will see later in the proof of  Theorem \ref{main_thm} that curvature plays a crucial role in our case. 

It should not be surprising to experts that the integrability exponent $(k^2+k+2)/2$ in Theorem \ref{main_thm} is the same as the sharp Fourier restriction exponent for the moment curve obtained by Drury \cite{Dr85}. 

Our argument works in a greater generality, for instance on $\Z/(N\Z)$; we will not expand the discussion here. \\

Finally we remark that the Agranovsky and Narayanan (\cite{AN04})  result was extended to fractal sets by Senthil Raani in \cite{SR14}. It is interesting to note that in this context, the exponent $\frac{2k}{s}$ is always essentially sharp, and $k$ need not be an integer. Roughly speaking, a random Cantor-type construction yields the sharp exponent, up to the endpoint. The key point is that given any $\alpha \in (0,d)$, it is possible to construct a compact random subset $E_{\alpha}$ of ${\mathbb R}^k$, of Minkowski (and Hausdorff) dimension $\alpha$, equipped with the Borel measure $\mu_{\alpha}$, such that for any $\epsilon>0$, there exists $C_{\epsilon}$ such that 
$$ |\widehat{\mu}_{\alpha}(\xi)| \leq C{(1+|\xi|)}^{-\alpha+\epsilon}.$$ 
Such sets are called Salem sets in honor of Raphael Salem who discovered them. On the other hand, it is not difficult to see that if $k<\frac{d}{2}$, no smooth $k$-dimensional submanifold of ${\mathbb R}^k$, $k \ge 2$, is a Salem set. One way to see this (see e.g. \cite{I99} and \cite{IL00}) is to show that such a decay rate would imply the $L^2$-restriction with the exponent which can be ruled out using the classical Knapp homogeneity argument. \\

\noindent {\bf Acknowledgement.} A. I. is supported in part by NSF DMS-2154232.
S. G. is partly supported by NSF-2044828. R. Z. is supported by NSF DMS-2207281(transferred from DMS-1856541), NSF DMS-2143989 and the Sloan Research Fellowship

\section{Proof of Theorem \ref{main_thm}}

\vskip.125in 

%For $\epsilon>0$, we define $\chi_{\epsilon}(t):=\chi(t/\epsilon)$. 
The proof will be carried out for the moment curve, but the argument easily carries over to a small perturbation of this case. For $t\in \R$, let $(\vector{e_1}(t), \dots, \vector{e_k}(t))$ denote the Frenet coordinate along the moment curve at the point $\vector{\gamma}(t)$. Here and below, all vectors are column vectors. Define a matrix 
\begin{equation}
M_t:=[\vector{e_1}(t), \dots, \vector{e_k}(t)].
\end{equation}
Define 
\begin{equation}
    D_{\epsilon}:=\diag[\epsilon, \dots, \epsilon^k],
\end{equation}
which refers to the diagonal matrix consisting of $\epsilon, \dots, \epsilon^k$ as diagonal elements. Let $\varphi: \R^k\to \R$ be a non-negative smooth bump function supported on $[-2\cdot 10^{-2}, 2\cdot 10^{-2}]^k$ which equals $1$ on $[-10^{-2}, 10^{-2}]^k$. Moreover, let $\widetilde{\varphi}: \R^k\to \R$ be a non-negative smooth bump function which equals $1$ in the $10^{-4}$-neighborhood of \eqref{moment_curve}.

Define an an-isotropic neighborhood of the moment curve by 
\begin{equation}
\Gamma_{\epsilon, t}:= \{\vector{\gamma}(t)+\epsilon_1 \vector{e_1}(t)+\dots \epsilon_k \vector{e_k}(t): |\epsilon_{k'}|\le \epsilon^{k'}, 1\le k'\le k\}
\end{equation}
Moreover, define 
\begin{equation}
\Gamma_{\epsilon}:=\bigcup_{t\in [0, 1]} \Gamma_{\epsilon, t}.
\end{equation}
Let $C_0$ be a large constant to be chosen. Define 
\begin{equation}
    \chi_{\epsilon, \iota}(\xi):=\varphi\Big(D_{\epsilon}^{-1}\cdot M^{-1}_{\iota \epsilon/C_0}\cdot (\xi-\vector{\gamma}(\iota \epsilon/C_0))\Big),
\end{equation}
and 
\begin{equation}
    \widetilde{\chi}_{\epsilon, \iota}(\xi):=\widetilde{\varphi}\Big(D_{\epsilon}^{-1}\cdot M^{-1}_{\iota \epsilon/C_0}\cdot (\xi-\vector{\gamma}
    (\iota\epsilon/C_0)
    )\Big),
\end{equation}
for an integer $\iota$. It is easy to see that if $C_0$ is chosen large enough, then 
\begin{equation}
    \sum_{\iota} \chi_{\epsilon, \iota} >0,
\end{equation}
at every point in $\Gamma_{\epsilon/C_0}$. Moreover, for every $\xi$ in 
\begin{equation}
\supp(\chi_{\epsilon, \iota})\cap \Gamma_{\epsilon/C_0},
\end{equation}
we have 
\begin{equation}
    \chi_{\epsilon, \iota}(\xi)\widetilde{\chi}_{\epsilon, \iota}(\xi)=\chi_{\epsilon, \iota}(\xi).
\end{equation}
Define 
\begin{equation}
    \eta_{\epsilon, \iota}(\xi):=\chi_{\epsilon, \iota}(\xi)\widetilde{\chi}_{\epsilon, \iota}(\xi) (\sum_{\iota} \chi_{\epsilon, \iota}(\xi))^{-1}. 
\end{equation}
So far we have obtained a partition of unity of $\Gamma_{\epsilon/C_0}$. In particular, 
\begin{equation}
    \eta_{\Gamma_ {\epsilon}}:=\sum_{\iota} \eta_{\epsilon, \iota} \equiv 1
\end{equation} 
on $\Gamma_{\epsilon/C_0}$. \\

% We pick a sequence of non-negative smooth functions $\{\chi_{\epsilon, \iota}\}_{\iota\in \Z}$ such that $\chi_{\epsilon, \iota}$ is supported in $\Gamma_{10\epsilon, \iota \epsilon}$ and that 
% \begin{equation}
% \chi_{\Gamma_{\epsilon}}:=\sum_{\iota} \chi_{\epsilon, \iota}\equiv 1 \text{ on } \Gamma_{\epsilon}. 
% \end{equation}
By the assumption that $\widehat{f}$ is supported on the moment curve, we obtain 
\begin{equation}
f=f* \widecheck{\eta}_{\Gamma_{\epsilon}},
\end{equation}
for every $\epsilon>0$. We claim that for every $p<(k^2+k+2)/2$, it holds that 
\begin{equation}
\|\widecheck{\eta}_{\Gamma_{\epsilon}}\|_{L^{p'}}\to 0 \text{ as } \epsilon\to 0.
\end{equation}
This will imply $\|f\|_{\infty}=0$ and therefore $f\equiv 0$. \\

To show the above claim, we just need to observe that 
\begin{equation}
\|\widecheck{\eta}_{\Gamma_{\epsilon}}\|_{L^{2}}\simeq \epsilon^{\frac{k(k+1)}{4}-\frac{1}{2}}. 
\end{equation}
Moreover, 
\begin{equation}
\|\widecheck{\eta}_{\Gamma_{\epsilon}}\|_{L^{1}}\le \sum_{\iota} \|\widecheck{\eta}_{\epsilon, \iota}\|_{L^{1}} \lesim \epsilon^{-1}. 
\end{equation}
Interpolation gives the desired result.

\bigskip 

\noindent Shaoming Guo: University of Wisconsin-Madison, shaomingguo@math.wisc.edu\\
Alex Iosevich: University of Rochester, iosevich@gmail.com\\
Ruixiang Zhang: UC Berkeley, ruixiang@berkeley.edu\\
Pavel Zorin-Kranich: University of Bonn, pzorinkranich@gmail.com

\end{document}